\documentclass[11pt,letterpaper,colorlinks=true,linkcolor=blue,citecolor=blue,urlcolor=blue]{scrartcl}
\usepackage[T1]{fontenc}
\usepackage[utf8]{inputenc}
\usepackage{lmodern}
\usepackage[english]{babel}
\usepackage[english]{isodate}
\usepackage{amsthm}
\usepackage{amsmath}
\usepackage{amssymb}
\usepackage{amsfonts}
\usepackage{mathtools}
\usepackage{mathrsfs}
\usepackage{dsfont}
\usepackage{upgreek}
\usepackage{xpatch}
\usepackage{bm}
\usepackage{bbm}
\usepackage{enumitem}
\usepackage{float}
\usepackage{footnotebackref}
\usepackage[totalwidth=480pt, totalheight=680pt]{geometry}
\usepackage{hyperref}

\theoremstyle{definition}
\newtheorem{definition}{Definition}[section]

\theoremstyle{plain}
\newtheorem{theorem}[definition]{Theorem}

\newtheorem{lemma}[definition]{Lemma}
\newtheorem{corollary}[definition]{Corollary}

\addtokomafont{disposition}{\normalfont\scshape}

\makeatletter
\xpatchcmd{\@sec@pppage}{
\bfseries}{
\normalfont\scshape\Large}{}{}
\makeatother

\setkomafont{section}{\Large}
\setkomafont{subsection}{\large}

\numberwithin{equation}{section}

\newcommand{\mybar}[1]{\makebox[0pt]{$\phantom{#1}\overline{\phantom{#1}}$}#1}

\begin{document}


\title{
\LARGE Diffusion processes as Wasserstein gradient flows via stochastic control of the volatility matrix\thanks{
I am grateful to Julio Backhoff-Veraguas and Walter Schachermayer for their useful comments. I also acknowledge support by the Austrian Science Fund (FWF) under grant P35519.
}}


\author{
\large Bertram Tschiderer\thanks{
Faculty of Mathematics, University of Vienna (email: \href{mailto: bertram.tschiderer@univie.ac.at}{bertram.tschiderer@univie.ac.at}).
} 
}


\date{}


\maketitle


\begin{abstract} \small \noindent \textsc{Abstract.} We consider a class of time-homogeneous diffusion processes on $\mathds{R}^{n}$ with common invariant measure but varying volatility matrices. In Euclidean space, we show via stochastic control of the diffusion coefficient that the corresponding flow of time-marginal distributions admits an entropic gradient flow formulation in the quadratic Wasserstein space if the volatility matrix of the diffusion is the identity. After equipping $\mathds{R}^{n}$ with a Riemannian metric, we prove that the diffusion process can be viewed as a gradient flow in the inherited Wasserstein space if the volatility matrix is the inverse of the underlying metric tensor. In the Euclidean case, our probabilistic result corresponds to the gradient flow characterization of the Fokker--Planck equation, first discovered in a seminal paper by Jordan, Kinderlehrer, and Otto. In the Riemannian setting, the corresponding result on the level of partial differential equations was established by Lisini, building on the metric theory developed by Ambrosio, Gigli, and Savar\'{e}.

\bigskip

\small \noindent \href{https://mathscinet.ams.org/mathscinet/msc/msc2020.html}{\textit{MSC 2020 subject classifications:}} Primary 60H30, 60G44; secondary 60J60, 94A17

\bigskip

\end{abstract}

\section{Introduction}

On a probability space $(\Omega,\mathcal{F},\mathds{P})$, equipped with filtration $\mathds{F} = (\mathcal{F}_{t})_{0 \leqslant t \leqslant T}$, we consider a class of time-homogeneous diffusion processes $(X_{t}^{\Sigma})_{0 \leqslant t \leqslant T}$ on $\mathds{R}^{n}$ of the form 
\begin{equation} \label{sde.sigma.v.lisini}
\textnormal{d}X_{t}^{\Sigma} 
= \Big( \operatorname{div} \Sigma(X_{t}^{\Sigma}) - \Sigma(X_{t}^{\Sigma}) \, \nabla V(X_{t}^{\Sigma}) \Big) \, \textnormal{d}t 
+ \sqrt{2\Sigma(X_{t}^{\Sigma})} \, \textnormal{d}B_{t},
\end{equation}
for $0 \leqslant t \leqslant T < \infty$; with initial condition $X_{0}^{\Sigma} = Y$, a square-integrable random vector with continuous probability density function $p_{0}^{\Sigma} \colon \mathds{R}^{n} \rightarrow (0,\infty)$, which is independent of the standard $n$-dimensional Brownian motion $(B_{t})_{0 \leqslant t \leqslant T}$. We denote by $P_{t}^{\Sigma} \coloneqq (X_{t}^{\Sigma})_{\#}\mathds{P}$ the distribution of the random vector $X_{t}^{\Sigma}$ and by $p_{t}^{\Sigma} \colon \mathds{R}^{n} \rightarrow (0,\infty)$ the corresponding probability density function. The functions $V \colon \mathds{R}^{n} \rightarrow \mathds{R}$ and $\Sigma \colon \mathds{R}^{n} \rightarrow \mathds{R}^{n \times n}$ are smooth and sufficiently well-behaved to guarantee the existence of a pathwise unique, strong solution of the stochastic differential equation \hyperref[sde.sigma.v.lisini]{(\ref*{sde.sigma.v.lisini})} in $L^{2}(\mathds{P})$, as well as a classical solution of the corresponding Fokker--Planck equation \cite{Fri75,Gar09,Ris96,Sch80} 
\begin{equation} \label{pde.sigma.v.lisini}
\partial_{t} p_{t}^{\Sigma}(x) 
= \operatorname{div}\Big(  \Sigma(x) \, \big( \nabla \log p_{t}^{\Sigma}(x) + \nabla V(x)\big) \, p_{t}^{\Sigma}(x) \Big),
\qquad (t,x) \in  (0,T) \times \mathds{R}^{n},
\end{equation}
which is also known as Kolmogorov forward equation \cite{Kol31}.

\smallskip

The function $\Sigma \colon \mathds{R}^{n} \rightarrow \mathds{R}^{n \times n}$ takes values in the set $S_{++}^{n}$ of symmetric and positive definite matrices and we denote by $\sqrt{\Sigma(x)} \in S_{++}^{n}$ the unique square root of $\Sigma(x) \in S_{++}^{n}$. The divergence of the matrix-valued function $\Sigma$ is the vector field $\operatorname{div} \Sigma \colon \mathds{R}^{n} \rightarrow \mathds{R}^{n}$ with components $(\operatorname{div} \Sigma)_{i} = \sum_{j=1}^{n} \partial_{j} \Sigma_{ij}$.

\smallskip

The aim of this paper is to study entropy dissipation and gradient flow properties of the diffusion equation \hyperref[sde.sigma.v.lisini]{(\ref*{sde.sigma.v.lisini})}. The novelty of our stochastic analysis approach is that it is based on stochastic control of the volatility matrix. Our primary inspiration derives from \cite{KST22}, where a similar endeavor has been established in the case when the volatility matrix $\Sigma$ is equal to the identity matrix.

\subsection{Invariant measure}

The stationary version 
\begin{equation} \label{pde.sigma.v.lisini.sv}
\operatorname{div}\Big(  \Sigma(x) \, \big( \nabla \log q(x) + \nabla V(x)\big) \, q(x) \Big) = 0,
\qquad x \in \mathds{R}^{n}
\end{equation}
of the Fokker--Planck equation \hyperref[pde.sigma.v.lisini]{(\ref*{pde.sigma.v.lisini})} has a strictly positive solution $q(x) \coloneqq \mathrm{e}^{-V(x)}$. We denote by $\mathrm{Q}$ the $\sigma$-finite measure on the Borel sets $\mathscr{B}(\mathds{R}^{n})$ having the function $q$ as density with respect to Lebesgue measure. This measure $\mathrm{Q}$ is invariant for the diffusion process $(X_{t}^{\Sigma})_{0 \leqslant t \leqslant T}$. If 
\begin{equation} \label{eq.gibbs.v.lisini}
\mathrm{Q}(\mathds{R}^{n}) 
= \int_{\mathds{R}^{n}} \mathrm{e}^{- V(x)} \, \textnormal{d}x < \infty,
\end{equation}
the measure $\mathrm{Q}$ can be normalized to become a probability measure (the so-called \textit{Gibbs measure} from statistical mechanics). We emphasize that we do not impose the condition \hyperref[eq.gibbs.v.lisini]{(\ref*{eq.gibbs.v.lisini})} and admit the possibility that the total mass $\mathrm{Q}(\mathds{R}^{n})$ is infinite. This allows for the important case when the function $V$ vanishes identically and $\mathrm{Q}$ is equal to Lebesgue measure on $\mathds{R}^{n}$. Instead, we require the weaker integrability condition 
\begin{equation} \label{int.cont.v.lisini}
\mathds{E}_{\mathrm{Q}}\Big[ \big( x \mapsto \mathrm{e}^{- \vert x \vert^{2}}\big) \Big] 
= \int_{\mathds{R}^{n}} \mathrm{e}^{- \vert x \vert^{2} - V(x)} \, \textnormal{d}x < \infty.
\end{equation}
If, for example, $V$ is non-negative, this condition is automatically satisfied. Furthermore, we stress that no convexity assumptions are imposed on the function $V$.

\smallskip 

We remark that, if $\operatorname{div} \Sigma = \Sigma \, \nabla V$, the process $(X_{t}^{\Sigma})_{0 \leqslant t \leqslant T}$ is a martingale. In this case, if $n=1$, the invariant measure $\mathrm{Q}$ has density $q(x) = \frac{1}{\Sigma(x)}$.

\subsection{Relative entropy and free energy}

Given the flow of probability measures $(P_{t}^{\Sigma})_{0 \leqslant t \leqslant T}$ defined by the diffusion process of \hyperref[sde.sigma.v.lisini]{(\ref*{sde.sigma.v.lisini})}, our goal is to study the \textit{relative entropy function}
\begin{equation} \label{rel.ent.function.lis.sig}
[0,T] \ni t \longmapsto H( P_{t}^{\Sigma} \, \vert \, \mathrm{Q}).
\end{equation}
The \textit{relative entropy} or \textit{Kullback--Leibler divergence} of $P_{t}^{\Sigma} = \operatorname{Law}(X_{t}^{\Sigma})$ with respect to the invariant measure $\mathrm{Q}$ is equal to
\begin{equation} \label{re.kld.lisini}
H( P_{t}^{\Sigma} \, \vert \, \mathrm{Q} ) 
= \mathds{E}_{\mathds{P}}\big[ \log \ell_{t}^{\Sigma}(X_{t}^{\Sigma})\big] 
= \int_{\mathds{R}^{n}} \log \bigg(\frac{p_{t}^{\Sigma}(x)}{q(x)}\bigg) \, p_{t}^{\Sigma}(x) \, \textnormal{d}x,
\end{equation}
where the \textit{likelihood ratio function} $(t,x) \mapsto \ell_{t}^{\Sigma}(x)$ is given by
\begin{equation} \label{eq.lis.lhrf.a}
\ell_{t}^{\Sigma}(x) \coloneqq \frac{p_{t}^{\Sigma}(x)}{q(x)} = p_{t}^{\Sigma}(x) \, \mathrm{e}^{V(x)}, 
\qquad (t,x) \in [0,T] \times \mathds{R}^{n}.
\end{equation}

\smallskip

The relative entropy \hyperref[re.kld.lisini]{(\ref*{re.kld.lisini})} can also be interpreted as a \textit{free energy}. Defining the \textit{free energy functional} as the sum of the internal and potential energy functionals via
\begin{equation} \label{def.fre.en.func.lisi}
p \longmapsto \mathscr{F}(p) \coloneqq \int_{\mathds{R}^{n}} p(x) \log p(x) \, \textnormal{d}x + \int_{\mathds{R}^{n}} V(x) \, p(x) \, \textnormal{d}x,
\end{equation}
for probability densities $p(\, \cdot \,)$ on $\mathds{R}^{n}$, we have the relation $H( P_{t}^{\Sigma} \, \vert \, \mathrm{Q} ) = \mathscr{F}(p_{t}^{\Sigma})$ between entropy and energy.

\smallskip

Since the diffusion process $(X_{t}^{\Sigma})_{0 \leqslant t \leqslant T}$ lies in $L^{2}(\mathds{P})$, the corresponding time-marginals $(P_{t}^{\Sigma})_{0 \leqslant t \leqslant T}$ are elements of $\mathscr{P}_{2}(\mathds{R}^{n})$, the collection of probability measures on $\mathscr{B}(\mathds{R}^{n})$ with finite second moments. Together with the integrability condition \hyperref[int.cont.v.lisini]{(\ref*{int.cont.v.lisini})}, this ensures that the relative entropy \hyperref[re.kld.lisini]{(\ref*{re.kld.lisini})} with respect to the $\sigma$-finite reference measure $\mathrm{Q}$ is well-defined and takes values in the interval $(-\infty,\infty]$; we refer to \cite[Appendix C]{KST20} or \cite[Section 3]{Leo14b} for the details. 

\subsection{Riemannian structure} \label{sec.1.4.lis.rs.asg}

We fix a function $A \in \mathcal{S}$, which then induces a smooth and symmetric metric tensor 
\[
G \colon \mathds{R}^{n} \longrightarrow \mathds{R}^{n \times n} \colon x \longmapsto G(x) \coloneqq A^{-1}(x).
\]
This metric tensor $G$ defines a Riemannian inner product $\langle \, \cdot \, ,  \cdot \, \rangle_{G}$ on $\mathds{R}^{n}$, given by 
\[
\mathds{R}^{n} \times \mathds{R}^{n} \longrightarrow \mathds{R} \colon 
(v,w) \longmapsto \langle v,w \rangle_{G}(x) \coloneqq  \langle G(x)v,w\rangle
\]
for $x \in \mathds{R}^{n}$. We denote the induced norm by $\vert  \cdot \vert_{G} \coloneqq \sqrt{\langle \, \cdot \, , \, \cdot \, \rangle_{G}}$. Equipped with this Riemannian structure, we can regard the pair $(\mathds{R}^{n},G)$ as a Riemannian manifold; we refer to \cite{Lee18} for an introduction to Riemannian manifolds. The metric tensor $G$ naturally induces a Riemannian metric $d_{G}$ on $(\mathds{R}^{n},G)$ defined by
\[
d_{G}(x,y) \coloneqq 
\inf_{\gamma \in \Gamma(x,y)} \bigg\{ 
\int_{0}^{1}  \sqrt{\Big\langle G\big(\gamma(t)\big) \dot{\gamma}(t) \, , \dot{\gamma}(t)\Big\rangle} \, \mathrm{d}t \bigg\}, \qquad x,y \in \mathds{R}^{n},
\]
where $\Gamma(x,y)$ denotes the collection of absolutely continuous curves $\gamma \colon [0,1] \rightarrow \mathds{R}^{n}$ with starting point $\gamma(0) = x$ and endpoint $\gamma(1) = y$. We write $\mathds{R}_{G}^{n}$ for the metric space $\mathds{R}^{n}$ equipped with this metric $d_{G}$. As a consequence of \hyperref[sigma.uec.lis]{(\ref*{sigma.uec.lis})}, $G$ satisfies a uniform ellipticity condition of the form
\begin{equation} \label{sigma.a.g}
\forall x, \xi \in \mathds{R}^{n} \colon \quad
C_{A}^{-1} \vert \xi \vert^{2} \leqslant \langle G(x) \xi , \xi \rangle \leqslant c_{A}^{-1} \vert \xi \vert^{2}.
\end{equation}
Note that $\mathscr{P}_{2}(\mathds{R}_{G}^{n}) = \mathscr{P}_{2}(\mathds{R}^{n})$ on account of \hyperref[sigma.a.g]{(\ref*{sigma.a.g})}. In other words, for a probability measure on $\mathds{R}^{n}$, it is immaterial whether the finite second moment property is tested with respect to the metric $d_{G}$ or the Euclidean metric, and by equivalence of norms the topology and the Borel sets remain the same.

\subsection{Relative Fisher information}

For a probability measure $P \in \mathscr{P}(\mathds{R}^{n})$ and vector fields $\boldsymbol{v},\boldsymbol{w} \colon \mathds{R}^{n} \rightarrow \mathds{R}^{n}$ we define the scalar product
\[
\langle \boldsymbol{v}, \boldsymbol{w} \rangle_{L_{G}^{2}(P)} 
\coloneqq \mathds{E}_{P} \big[ \langle \boldsymbol{v}, \boldsymbol{w} \rangle_{G} \big] 
= \int_{\mathds{R}^{n}} \langle G\boldsymbol{v}, \boldsymbol{w} \rangle \, \textnormal{d}P
\]
and the induced norm $\Vert \boldsymbol{v} \Vert_{L_{G}^{2}(P)} \coloneqq \sqrt{\langle \boldsymbol{v}, \boldsymbol{v} \rangle_{L_{G}^{2}(P)}}$.

\smallskip

The gradient of a smooth function $f \colon \mathds{R}^{n} \rightarrow \mathds{R}$ on the Riemannian manifold $(\mathds{R}^{n},G)$ is given by $\nabla_{G}f \coloneqq A \nabla f$ (see, e.g., \cite{Lee18}).

\smallskip

In our Riemannian setting, the \textit{relative Fisher information} (see, e.g., \cite{OV00,Vil03}) of the probability measure $P_{t}^{\Sigma} = \operatorname{Law}(X_{t}^{\Sigma})$ with respect to the invariant measure $\mathrm{Q}$ is equal to
\[
I_{G}( P_{t}^{\Sigma} \, \vert \, \mathrm{Q}) 
\coloneqq \Vert \nabla_{G} \log \ell_{t}^{\Sigma} \Vert_{L_{G}^{2}(P_{t}^{\Sigma})}^{2}
= \int_{\mathds{R}^{n}} \big\langle \nabla \log \ell_{t}^{\Sigma}  ,  A \, \nabla \log \ell_{t}^{\Sigma} \big\rangle(x)  \, p_{t}^{\Sigma}(x) \, \textnormal{d}x.
\]
We can also express the relative Fisher information as a $\mathds{P}$-expectation, namely
\begin{equation} \label{rel.fis.inf.det.pexp.sig}
I_{G}( P_{t}^{\Sigma} \, \vert \, \mathrm{Q}) 
= \mathds{E}_{\mathds{P}}\big[ \vert \nabla_{G} \log \ell_{t}^{\Sigma}(X_{t}^{\Sigma})\vert_{G}^{2}  \big].
\end{equation}
For the classical definition of Fisher information in the information theory literature we refer to the books \cite{CL98, CT06}.

\subsection{Wasserstein distance}

The quadratic Wasserstein distance (see, e.g., \cite{AGS08, Vil03}) between two probability measures $P_{1}$ and $P_{2}$ in $\mathscr{P}_{2}(\mathds{R}^{n})$, with $\mathds{R}^{n}$ being equipped with the Riemannian metric $d_{G}$, is defined by 
\begin{equation} \label{def.eq.wtwog.lis.was}
W_{2,G}(P_{1},P_{2}) \coloneqq 
\sqrt{\inf_{\pi \in \Pi(P_{1},P_{2})} \int_{\mathds{R}^{n} \times \mathds{R}^{n}} d_{G}(x,y)^{2} \, \textnormal{d} \pi(x,y)} \, ,
\end{equation}
where $\Pi(P_{1},P_{2})$ denotes the set of all couplings $\pi \in \mathscr{P}(\mathds{R}^{n} \times \mathds{R}^{n})$ between $P_{1}$ and $P_{2}$, i.e., probability measures $\pi$ on $\mathds{R}^{n} \times \mathds{R}^{n}$ with first marginal $P_{1}$ and second marginal $P_{2}$. 

\subsection{Class of admissible volatilities}

\begin{definition} \label{def.cl.of.ad.vol.sig.lis} We define the \textit{class of admissible volatilities} $\mathcal{S}$ as the collection of all matrix-valued functions $\Sigma \colon \mathds{R}^{n} \rightarrow \mathds{R}^{n \times n}$ such that
\begin{enumerate}[label=(\roman*)] 
\item $\Sigma(x)$ is a symmetric and positive definite matrix, for every $x \in \mathds{R}^{n}$,
\item the component functions $\mathds{R}^{n} \ni x \mapsto \Sigma_{ij}(x) \in \mathds{R}$ are smooth, for all $i,j \in \{1, \ldots, n\}$,
\item for some constants $c,C > 0$ the uniform ellipticity condition 
\begin{equation} \label{sigma.uec.lis}
\forall x, \xi \in \mathds{R}^{n} \colon \quad
c \vert \xi \vert^{2} \leqslant \langle \Sigma(x) \xi , \xi \rangle \leqslant C \vert \xi \vert^{2} 
\end{equation}
is satisfied,
\item the stochastic differential equation \hyperref[sde.sigma.v.lisini]{(\ref*{sde.sigma.v.lisini})} admits a pathwise unique, strong solution $(X_{t}^{\Sigma})_{0 \leqslant t \leqslant T}$, which is bounded in $L^{2}(\mathds{P})$,
\item the Fokker--Planck equation \hyperref[pde.sigma.v.lisini]{(\ref*{pde.sigma.v.lisini})} has a classical solution $(p_{t}^{\Sigma})_{0 \leqslant t \leqslant T}$,
\item the curve of probability measures $(P_{t}^{\Sigma})_{0 \leqslant t \leqslant T}$ given by $P_{t}^{\Sigma} = \operatorname{Law}(X_{t}^{\Sigma})$ is absolutely continuous on the quadratic Wasserstein space $(\mathscr{P}_{2}(\mathds{R}^{n}),W_{2,G})$, i.e., there exists $m \in L^{1}([0,T])$ such that
\[
W_{2,G}(P_{t_{0}}^{\Sigma},P_{t}^{\Sigma}) \leqslant \int_{t_{0}}^{t} m(u) \, \textnormal{d}u
\]
for all $0 \leqslant t_{0} \leqslant t \leqslant T$,
\item \label{def.cl.of.ad.vol.sig.lis.vi} and the initial relative entropy $H( P_{0}^{\Sigma} \, \vert \, \mathrm{Q})$ is finite.
\end{enumerate}
\end{definition}

\subsection{Main results}

Under the assumptions of this section, we can state our first main result.

\begin{theorem} \label{main.one.lis.det} For every $\Sigma \in \mathcal{S}$ and $0 \leqslant t_{0} \leqslant t \leqslant T$ we have
\begin{equation} \label{main.one.lis.det.01}
\big\vert H( P_{t}^{\Sigma} \, \vert \, \mathrm{Q}) - H( P_{t_{0}}^{\Sigma} \, \vert \, \mathrm{Q})\big\vert
\leqslant \int_{t_{0}}^{t} \sqrt{I_{G}( P_{u}^{\Sigma} \, \vert \, \mathrm{Q})} \ \lim_{h \rightarrow 0} \frac{W_{2,G}(P_{u+h}^{\Sigma},P_{u}^{\Sigma})}{\vert h \vert} \, \textnormal{d}u,
\end{equation}
with equality if the diffusion matrix $\Sigma$ is a constant multiple of $A$. In the case $\Sigma = A$ we have
\begin{equation} \label{main.one.lis.det.02}
H( P_{t_{0}}^{A} \, \vert \, \mathrm{Q}) -  H( P_{t}^{A} \, \vert \, \mathrm{Q}) 
= \tfrac{1}{2} \int_{t_{0}}^{t} I_{G}( P_{u}^{A} \, \vert \, \mathrm{Q}) \, \textnormal{d}u
+ \tfrac{1}{2} \int_{t_{0}}^{t} \bigg( \lim_{h \rightarrow 0} \frac{W_{2,G}(P_{u+h}^{A},P_{u}^{A})}{\vert h \vert} \bigg)^{2} \, \textnormal{d}u.
\end{equation}
\end{theorem}

The message of \hyperref[main.one.lis.det]{Theorem \ref*{main.one.lis.det}} is that the ``flow of probability measures'' $(P_{t}^{A})_{0 \leqslant t \leqslant T}$, defined by the stochastic differential equation
\begin{equation} \label{sde.sigma.v.lisini.A}
\textnormal{d}X_{t}^{A} 
= \Big( \operatorname{div} A(X_{t}^{A}) - A(X_{t}^{A}) \, \nabla V(X_{t}^{A}) \Big) \, \textnormal{d}t 
+ \sqrt{2A(X_{t}^{A})} \, \textnormal{d}B_{t},
\end{equation}
is the \textit{gradient flow} of the \textit{relative entropy functional}
\begin{equation} \label{eq.rel.ent.func.lis.}
\mathscr{P}_{2}(\mathds{R}^{n}) \ni P \longmapsto H( P \, \vert \, \mathrm{Q})
\end{equation}
with respect to the quadratic Wasserstein distance $W_{2,G}$. The gradient flow property in this Riemannian setting was established by Lisini in \cite{Lis09}, where nonlinear diffusion equations with variable coefficients are identified as gradient flows in Wasserstein spaces with respect to a Riemannian metric. The approach of \cite{Lis09} is based on the theory of minimizing movements and the theory of curves of maximal slope in metric spaces, taking as a reference \cite{AGS08}. While in \cite{Lis09} it is crucial to study the Fokker--Planck equation 
\begin{equation} \label{pde.sigma.v.lisini.A.a}
\partial_{t} p_{t}^{A}(x) 
= \operatorname{div}\Big(  A(x) \, \big( \nabla \log p_{t}^{A}(x) + \nabla V(x)\big) \, p_{t}^{A}(x) \Big),
\end{equation}
our approach is based on the corresponding stochastic differential equation \hyperref[sde.sigma.v.lisini.A]{(\ref*{sde.sigma.v.lisini.A})} via stochastic control of the diffusion coefficient. More precisely, we fix $A \in \mathcal{S}$, and for every $\Sigma$ in the class of admissible volatilities $\mathcal{S}$, we vary the diffusion coefficient by considering a stochastic differential equation of the form \hyperref[sde.sigma.v.lisini]{(\ref*{sde.sigma.v.lisini})}. After some stochastic analysis and taking expectations this yields the inequality \hyperref[main.one.lis.det.01]{(\ref*{main.one.lis.det.01})}. In the language of analysis in metric spaces, the inequality \hyperref[main.one.lis.det.01]{(\ref*{main.one.lis.det.01})} has the interpretation that the square root of the \textit{relative Fisher information functional}
\begin{equation} \label{st.up.gr.re.fi.g.lis}
P \longmapsto \sqrt{I_{G}( P \, \vert \, \mathrm{Q})}
= \Vert \nabla_{G} \log (\tfrac{p}{q}) \Vert_{L_{G}^{2}(P)}
= \sqrt{\int_{\mathds{R}^{n}} \big\langle \nabla \log (\tfrac{p}{q})  ,  A \, \nabla \log (\tfrac{p}{q}) \big\rangle(x)  \, p(x) \, \textnormal{d}x}
\end{equation}
is a \textit{strong upper gradient} (see \cite[Definition 1.2.1]{AGS08}, \cite[Definition 3.2]{Lis09}) for the relative entropy functional \hyperref[eq.rel.ent.func.lis.]{(\ref*{eq.rel.ent.func.lis.})}. The domain of these two functionals consists of the collection of absolutely continuous probability measures $P \in \mathscr{P}_{2}(\mathds{R}^{n})$, having a probability density function $p(\, \cdot \,)$ with respect to $n$-dimensional Lebesgue measure. In the case $\Sigma = A$, we obtain the so-called \textit{energy identity} \hyperref[main.one.lis.det.02]{(\ref*{main.one.lis.det.02})} (see \cite[Remark 1.3.3]{AGS08}, \cite[Definition 3.3]{Lis09}). This shows that $(P_{t}^{A})_{0 \leqslant t \leqslant T}$ is a \textit{curve of maximal slope} (see \cite[Definition 1.3.2]{AGS08}, \cite[Definition 3.3]{Lis09}) for the relative entropy functional \hyperref[eq.rel.ent.func.lis.]{(\ref*{eq.rel.ent.func.lis.})} with respect to its strong upper gradient \hyperref[st.up.gr.re.fi.g.lis]{(\ref*{st.up.gr.re.fi.g.lis})} in the Wasserstein space $(\mathscr{P}_{2}(\mathds{R}^{n}),W_{2,G})$. According to \cite[Theorem 11.1.3]{AGS08}, \textit{curves of maximal slope coincide with gradient flows}, which verifies the stated gradient flow property of $(P_{t}^{A})_{0 \leqslant t \leqslant T}$.

\medskip

We leave now the Riemannian setting and return to the standard Euclidean structure on $\mathds{R}^{n}$. In other words, we fix the identity matrix $A = I_{n} \in \mathcal{S}$ and consider the quadratic Wasserstein distance
\begin{equation} \label{eq.quad.was.dist.lis}
W_{2}(P_{1},P_{2}) =
\sqrt{\inf_{\pi \in \Pi(P_{1},P_{2})} \int_{\mathds{R}^{n} \times \mathds{R}^{n}} \vert x-y \vert^{2} \, \textnormal{d} \pi(x,y)} \, , \qquad P_{1}, P_{2} \in \mathscr{P}_{2}(\mathds{R}^{n}),
\end{equation}
together with the relative Fisher information functional
\[
P \longmapsto I( P \, \vert \, \mathrm{Q}) 
= \Vert \nabla \log (\tfrac{p}{q}) \Vert_{L^{2}(P)}^{2}
= \int_{\mathds{R}^{n}} \vert \nabla \log (\tfrac{p(x)}{q(x)}) \vert^{2}  \, p(x) \, \textnormal{d}x.
\]
As an immediate consequence of \hyperref[main.one.lis.det]{Theorem \ref*{main.one.lis.det}}, we obtain the following result.

\begin{corollary} \label{main.one.lis.det.cor} For every $\Sigma \in \mathcal{S}$ and $0 \leqslant t_{0} \leqslant t \leqslant T$ we have
\[
\big\vert H( P_{t}^{\Sigma} \, \vert \, \mathrm{Q}) - H( P_{t_{0}}^{\Sigma} \, \vert \, \mathrm{Q})\big\vert
\leqslant \int_{t_{0}}^{t} \sqrt{I( P_{u}^{\Sigma} \, \vert \, \mathrm{Q})} \ \lim_{h \rightarrow 0} \frac{W_{2}(P_{u+h}^{\Sigma},P_{u}^{\Sigma})}{\vert h \vert} \, \textnormal{d}u,
\]
with equality if the diffusion matrix $\Sigma$ is a constant multiple of the identity matrix. In the case $\Sigma = I_{n}$ we have 
\[
H( P_{t_{0}} \, \vert \, \mathrm{Q}) -  H( P_{t} \, \vert \, \mathrm{Q}) 
= \tfrac{1}{2} \int_{t_{0}}^{t} I( P_{u} \, \vert \, \mathrm{Q}) \, \textnormal{d}u
+ \tfrac{1}{2} \int_{t_{0}}^{t} \bigg( \lim_{h \rightarrow 0} \frac{W_{2}(P_{u+h},P_{u})}{\vert h \vert} \bigg)^{2} \, \textnormal{d}u.
\]
\end{corollary}

Here, the probability measures $(P_{t})_{0 \leqslant t \leqslant T}$ denote the time-marginal distributions of the stochastic differential equation of Langevin--Smoluchowski type
\begin{equation} \label{sde.identity.v.lisini}
\textnormal{d}X_{t}
= -  \nabla V(X_{t}) \, \textnormal{d}t + \sqrt{2} \, \textnormal{d}B_{t}.
\end{equation}
The corresponding probability density functions $(p_{t})_{0 \leqslant t \leqslant T}$ satisfy the Fokker--Planck equation
\begin{equation} \label{pde.identity.v.lisini}
\partial_{t} p_{t}(x) 
= \operatorname{div}\Big(   \big( \nabla \log p_{t}(x) + \nabla V(x)\big) \, p_{t}(x) \Big).
\end{equation}
Clearly the equations \hyperref[sde.identity.v.lisini]{(\ref*{sde.identity.v.lisini})} and \hyperref[pde.identity.v.lisini]{(\ref*{pde.identity.v.lisini})} are obtained by taking $A$ to be the $n$-dimensional identity matrix in \hyperref[sde.sigma.v.lisini.A]{(\ref*{sde.sigma.v.lisini.A})} and \hyperref[pde.sigma.v.lisini.A.a]{(\ref*{pde.sigma.v.lisini.A.a})}, respectively. By analogy with the Riemannian setting considered above, \hyperref[main.one.lis.det.cor]{Corollary \ref*{main.one.lis.det.cor}} has a gradient flow interpretation in Euclidean space. Namely, it establishes that the flow of probability measures $(P_{t})_{0 \leqslant t \leqslant T}$ is the gradient flow of the relative entropy functional \hyperref[eq.rel.ent.func.lis.]{(\ref*{eq.rel.ent.func.lis.})} with respect to the quadratic Wasserstein distance $W_{2}$ on the space $\mathscr{P}_{2}(\mathds{R}^{n})$. The gradient flow property in this context was first established in \cite{JKO98}.

\bigskip

The special feature of our approach is that it is mainly based on stochastic analysis. In fact, most of \hyperref[main.one.lis.det]{Theorem \ref*{main.one.lis.det}}, which is a deterministic result, will be a consequence of a stronger trajectorial result, \hyperref[main.two.lis.tra]{Theorem \ref*{main.two.lis.tra}}. This constitutes the second main result of this paper, it is stated and proved in \hyperref[dotrep.chap.lis]{Section \ref*{dotrep.chap.lis}}. Trajectorial approaches of a similar kind were initiated by \cite{FJ16, KST22}. Several important consequences and ramifications of \hyperref[main.two.lis.tra]{Theorem \ref*{main.two.lis.tra}} are developed in \hyperref[cottt.lis.ram.sig.a]{Subsection \ref*{cottt.lis.ram.sig.a}}.

\smallskip

In \hyperref[ch.lis.tpotgrp]{Section \ref*{ch.lis.tpotgrp}} we prove the gradient flow property as stated in \hyperref[main.one.lis.det]{Theorem \ref*{main.one.lis.det}}. Most of the proof will be a consequence of \hyperref[main.two.lis.tra]{Theorem \ref*{main.two.lis.tra}}. As a final ingredient, we need to compute the derivative of the Wasserstein distance 
\begin{equation} \label{eq.fu.wmd.lis.a}
f(u) \coloneqq \lim_{h \rightarrow 0} \frac{W_{2,G}(P_{u+h}^{\Sigma},P_{u}^{\Sigma})}{\vert h \vert},
\end{equation}
a quantity appearing in \hyperref[main.one.lis.det.01]{(\ref*{main.one.lis.det.01})} as well as \hyperref[main.one.lis.det.02]{(\ref*{main.one.lis.det.02})}. The derivative \hyperref[eq.fu.wmd.lis.a]{(\ref*{eq.fu.wmd.lis.a})} is known as a \textit{metric derivative} in \cite{AGS08, Lis09}. By relying on the metric theory of \cite{AGS08, Lis09} we show in \hyperref[sec.3.lis.dotwd]{Subsection \ref*{sec.3.lis.dotwd}} that
\[
f(u) = \sqrt{\mathds{E}_{\mathds{P}}\Big[ \big\vert \big(\Sigma \, G \, \nabla_{G} \log \ell_{u}^{\Sigma}\big)(X_{u}^{\Sigma})\big\vert_{G}^{2}  \Big]}.
\]
Finally, in \hyperref[section4.lis.tpot11]{Subsection \ref*{section4.lis.tpot11}} we merge our findings and provide the proof of \hyperref[main.one.lis.det]{Theorem \ref*{main.one.lis.det}}.

\subsection{Literature review}

In the seminal paper \cite{JKO98}, Jordan, Kinderlehrer, and Otto established a new relationship between a particular class of Fokker--Planck equations, the associated free energy functionals, and the quadratic Wasserstein distance. More precisely, they considered the flow of probability densities described by a stochastic differential equation of Langevin--Smoluchowski type as in \hyperref[sde.identity.v.lisini]{(\ref*{sde.identity.v.lisini})}, for which the drift term is given by the gradient of a potential and the diffusion coefficient is constant. The authors of \cite{JKO98} constructed a time discrete, iterative variational scheme, whose solutions converge to the solution of the corresponding Fokker--Planck equation \hyperref[pde.identity.v.lisini]{(\ref*{pde.identity.v.lisini})}. Associated with this Fokker--Planck equation is a free energy functional, modeled as the sum of internal and potential energy as in \hyperref[def.fre.en.func.lisi]{(\ref*{def.fre.en.func.lisi})}, and defined on the set of probability densities with finite second moments. The path-breaking result shown in \cite{JKO98} is that the slope of this free energy functional along the Fokker--Planck probability density flow is the steepest possible, if one chooses the quadratic Wasserstein distance \hyperref[eq.quad.was.dist.lis]{(\ref*{eq.quad.was.dist.lis})} as metric. In other words, this flow of probability densities can be regarded as a ``curve of steepest descent'', or a ``gradient flux'', for the free energy functional with respect to the Wasserstein distance. 

\smallskip

The point of view taken in \cite{JKO98} was further developed in the paper \cite{Ott01}, where Otto showed the Wasserstein gradient flow structure of the porous medium equation and introduced a formal Riemannian calculus on the space of probability density functions, which was later dubbed ``Otto calculus'' in Chapter 15 of \cite{Vil09}; see also Chapter 8 of \cite{Vil03}. As just one example, the results by Otto and Villani on transport inequalities in \cite{OV00} show the power of this heuristic method. For rigorous constructions we refer to the research papers \cite{CG03,CMV06,Lot08,LV09,Oht09,Stu06a,Stu06b} and the textbook \cite{AGS08} by Ambrosio, Gigli, and Savar\'{e}.

\smallskip

For the general and rigorous definition of gradient flows in metric spaces we again refer to the book \cite{AGS08}, which is based on the theory of minimizing movements and the theory of curves of maximal slope in metric spaces, building on previous work like \cite{DG93,DGMT80,DMT85}. For literature on optimal transport and gradient flows we also refer to \cite{Vil03, Vil09, AG13, San15}.

\smallskip

The stochastic analysis approach of this work stems from \cite{FJ16} and \cite{KST22}. Such trajectorial approaches, including applications of time-reversal, were also carried forward in \cite{KT22, KMS20, TY23, KY22, CL22, JK23}. For the connection of gradient flows with large deviations we mention \cite{ADPZ13, Fat16}; further recent related work includes \cite{BCGL20, HRSS21}.

\section{Dynamics of the relative entropy process} \label{dotrep.chap.lis}

In this section we present and prove our second main result, \hyperref[main.two.lis.tra]{Theorem \ref*{main.two.lis.tra}}. In contrast to \hyperref[main.one.lis.det]{Theorem \ref*{main.one.lis.det}} it is not of a deterministic but of a trajectorial nature. Most of \hyperref[main.one.lis.det]{Theorem \ref*{main.one.lis.det}} will then follow from \hyperref[main.two.lis.tra]{Theorem \ref*{main.two.lis.tra}} by taking expectations.

\smallskip

Our eventual goal is to study the time evolution of the relative entropy function \hyperref[rel.ent.function.lis.sig]{(\ref*{rel.ent.function.lis.sig})}. We will do this in a trajectorial fashion by modelling the relative entropy function as a stochastic process. This idea goes back to \cite{FJ16} and was utilized in the context of gradient flows in \cite{KST22}. In the latter work a stochastic analysis approach to the characterization of diffusions of Langevin--Smoluchowski type \hyperref[sde.identity.v.lisini]{(\ref*{sde.identity.v.lisini})} as entropic gradient flows in Wasserstein space was developed. The main difference to our current work is that \cite{KST22} is formulated in the setting of \cite{JKO98}, where the diffusion coefficient is constant and the underlying metric is Euclidean. Along the lines of \cite{Lis09}, we allow here for a variable diffusion coefficient as well as a Riemannian geometry (determined by $A \in \mathcal{S}$ as in \hyperref[sec.1.4.lis.rs.asg]{Subsection \ref*{sec.1.4.lis.rs.asg}}). As a consequence, while in \cite{KST22} the drift coefficient $\nabla V$ is controlled (or, put another way, perturbed), we vary the diffusion coefficient $\Sigma$ in the class of admissible volatilities $\mathcal{S}$. In particular, we provide probabilistic interpretations of the corresponding results in \cite{Lis09} by considering the associated diffusion process \hyperref[sde.sigma.v.lisini.A]{(\ref*{sde.sigma.v.lisini.A})} rather than the underlying Fokker--Planck equation \hyperref[pde.sigma.v.lisini.A.a]{(\ref*{pde.sigma.v.lisini.A.a})}. From this perspective, our work can be viewed as a generalization of \cite{KST22}, in a similar way as \cite{Lis09} generalizes the setting of \cite{JKO98} and some of the results in \cite{AGS08}. Clearly, by choosing the underlying metric to be Euclidean (i.e., by fixing the identity matrix $A = I_{n} \in \mathcal{S}$ as a reference), we recover the results of \cite{JKO98, KST22}; recall \hyperref[main.one.lis.det.cor]{Corollary \ref*{main.one.lis.det.cor}} in this context. We also remark that the setting of \cite{Lis09} allows for more general nonlinear partial differential equations, as opposed to the linear case \hyperref[pde.sigma.v.lisini.A.a]{(\ref*{pde.sigma.v.lisini.A.a})}.
 
\medskip

We fix $\Sigma \in \mathcal{S}$ and consider the diffusion process $(X_{t}^{\Sigma})_{0 \leqslant t \leqslant T}$, governed by the stochastic differential equation \hyperref[sde.sigma.v.lisini]{(\ref*{sde.sigma.v.lisini})}, with time-marginal distributions $P_{t}^{\Sigma} = \operatorname{Law}(X_{t}^{\Sigma})$. With the above motivation in mind, we consider the likelihood ratio function \hyperref[eq.lis.lhrf.a]{(\ref*{eq.lis.lhrf.a})} and its logarithm as stochastic processes. The resulting random objects are the \textit{likelihood ratio process} (or \textit{Radon--Nikod\'{y}m derivative process}) 
\begin{equation} \label{llhr.pr.itfdot.lis.s}
\frac{\textnormal{d}P_{t}^{\Sigma}}{\textnormal{d}\mathrm{Q}}(X_{t}^{\Sigma}) 
= \ell_{t}^{\Sigma}(X_{t}^{\Sigma}) 
= \frac{p_{t}^{\Sigma}(X_{t}^{\Sigma})}{q(X_{t}^{\Sigma})} 
= p_{t}^{\Sigma}(X_{t}^{\Sigma}) \, \mathrm{e}^{V(X_{t}^{\Sigma})}, 
\qquad 0 \leqslant t \leqslant T
\end{equation}
and the \textit{relative entropy process} 
\begin{equation} \label{eq.rel.ent.pro.fdot.lis}
\log \ell_{t}^{\Sigma}(X_{t}^{\Sigma}) 
= \log \bigg( \frac{p_{t}^{\Sigma}(X_{t}^{\Sigma})}{q(X_{t}^{\Sigma})} \bigg) 
= \log p_{t}^{\Sigma}(X_{t}^{\Sigma}) + V(X_{t}^{\Sigma}),
\qquad 0 \leqslant t \leqslant T.
\end{equation}
Since taking the $\mathds{P}$-expectation in \hyperref[eq.rel.ent.pro.fdot.lis]{(\ref*{eq.rel.ent.pro.fdot.lis})} leads to the relative entropy as in \hyperref[re.kld.lisini]{(\ref*{re.kld.lisini})}, the idea is to compute the dynamics of the relative entropy process. More precisely, we will determine the stochastic differential of the process \hyperref[eq.rel.ent.pro.fdot.lis]{(\ref*{eq.rel.ent.pro.fdot.lis})} under the probability measure $\mathds{P}$. Motivated from results like \cite{DPP89, Pav89}, and most notably the work \cite{FJ16} by Fontbona and Jourdain, as well as the stochastic approach of \cite{KST22}, we expect that these dynamics will be most transparent when regarded in the backward direction of time. Thus we consider the time-reversed process 
\begin{equation} \label{eq.trxsig.lis}
\mybar{X}_{s}^{\Sigma} 
\coloneqq X_{T-s}^{\Sigma}, 
\qquad 0 \leqslant s \leqslant T.
\end{equation}
Due to the strong solvability of the stochastic differential equation \hyperref[sde.sigma.v.lisini]{(\ref*{sde.sigma.v.lisini})}, we can assume that the filtration $\mathds{F} = (\mathcal{F}_{t})_{0 \leqslant t \leqslant T}$ is given by
\[
\mathcal{F}_{t} 
= \sigma\big( X_{u}^{\Sigma}, B_{u} \colon \, 0 \leqslant u \leqslant t \big)
\]
modulo $\mathds{P}$-augmentation. We define the \textit{backwards filtration} $\mybar{\mathds{F}} = (\mybar{\mathcal{F}}_{s})_{0 \leqslant s \leqslant T}$ by
\[
\mybar{\mathcal{F}}_{s} 
\coloneqq \sigma\big( \mybar{X}_{u}^{\Sigma}, B_{T-u} - B_{T-s} \colon \, 0 \leqslant u \leqslant s \big)
\]
modulo $\mathds{P}$-augmentation.

\smallskip

The following \hyperref[main.two.lis.tra]{Theorem \ref*{main.two.lis.tra}} describes the semimartingale dynamics of the time-reversed relative entropy process 
\begin{equation} \label{eq.smdottrrop.lis}
\log \mybar{\ell}_{s}^{\Sigma}(\mybar{X}_{s}^{\Sigma}) 
\coloneqq \log \ell_{T-s}^{\Sigma}(X_{T-s}^{\Sigma}), 
\qquad 0 \leqslant s \leqslant T
\end{equation}
in terms of the \textit{cumulative relative Fisher information process} 
\begin{equation} \label{cum.rel.fis.inf.pro.lis.sig}
\mybar{F}_{s}^{\Sigma} \coloneqq
\int_{0}^{s} \Big\langle \nabla \log \mybar{\ell}_{u}^{\Sigma} \, , \, \Sigma  \, \nabla \log \mybar{\ell}_{u}^{\Sigma} \Big\rangle(\mybar{X}_{u}^{\Sigma}) \, \textnormal{d}u,
\end{equation}
for $0 \leqslant s \leqslant T$. In particular, we will see that the difference between the processes of \hyperref[eq.smdottrrop.lis]{(\ref*{eq.smdottrrop.lis})} and \hyperref[cum.rel.fis.inf.pro.lis.sig]{(\ref*{cum.rel.fis.inf.pro.lis.sig})} is a martingale. The proof of the theorem below is given in \hyperref[sec.2.1.lis.potmt]{Subsection \ref*{sec.2.1.lis.potmt}}.

\begin{theorem} \label{main.two.lis.tra} Let $\Sigma \in \mathcal{S}$, fix $A \in \mathcal{S}$ and set $G = A^{-1}$. On the filtered probability space $(\Omega,\mybar{\mathcal{F}}_{T},\mybar{\mathds{F}},\mathds{P})$ we have
\begin{equation} \label{main.two.lis.tra.a}
\mathds{E}_{\mathds{P}}\big[ \, \mybar{F}_{T}^{\Sigma} \, \big] 
= \int_{0}^{T} \mathds{E}_{\mathds{P}}\Big[ \Big\langle \nabla \log \ell_{t}^{\Sigma} \, , \, \Sigma  \, \nabla \log \ell_{t}^{\Sigma} \Big\rangle(X_{t}^{\Sigma}) \Big] \, \textnormal{d}t
< \infty
\end{equation}
and the process
\begin{equation} \label{main.two.lis.tra.b}
\mybar{M}_{s}^{\Sigma} 
\coloneqq \big( \log \mybar{\ell}_{s}^{\Sigma} (\mybar{X}_{s}^{\Sigma} ) - \log \ell_{T}^{\Sigma} (X_{T}^{\Sigma} )  \big) - \mybar{F}_{s}^{\Sigma},
\qquad 0 \leqslant s \leqslant T
\end{equation}
is a martingale bounded in $L^{2}(\mathds{P})$, with representation as a stochastic integral
\begin{equation} \label{main.two.lis.tra.c}
\mybar{M}_{s}^{\Sigma}  
= \int_{0}^{s} \Big\langle \nabla \log \mybar{\ell}_{u}^{\Sigma}(\mybar{X}_{u}^{\Sigma} ) \, , \sqrt{2\Sigma(\mybar{X}_{u}^{\Sigma} )} \, \textnormal{d}\mybar{B}_{u} \Big\rangle, 
\end{equation}
where the process $(\mybar{B}_{s})_{0 \leqslant s \leqslant T}$ is $\mybar{\mathds{F}}$-Brownian motion \textnormal{(}see \textnormal{\hyperref[sec.2.1.lis.potmt.lem.bm]{Lemma \ref*{sec.2.1.lis.potmt.lem.bm}}}\textnormal{)}.
\end{theorem}

We now consider the important special case $\Sigma = A$. In this situation, the cumulative relative Fisher information process \hyperref[cum.rel.fis.inf.pro.lis.sig]{(\ref*{cum.rel.fis.inf.pro.lis.sig})} takes the form 
\[
\mybar{F}_{s}^{A} =
\int_{0}^{s} \Big\langle \nabla \log \mybar{\ell}_{u}^{A} \, , \, A  \, \nabla \log \mybar{\ell}_{u}^{A} \Big\rangle(\mybar{X}_{u}^{\Sigma}) \, \textnormal{d}u
= \int_{0}^{s} \vert \nabla_{G} \log \mybar{\ell}_{u}^{A}(\mybar{X}_{u}^{A})\vert_{G}^{2} \, \textnormal{d}u
\]
and in light of \hyperref[rel.fis.inf.det.pexp.sig]{(\ref*{rel.fis.inf.det.pexp.sig})} we have
\[
\mathds{E}_{\mathds{P}}\big[ \, \mybar{F}_{T}^{A} \, \big] 
= \int_{0}^{T} \mathds{E}_{\mathds{P}}\Big[ \vert \nabla_{G} \log \ell_{t}^{A}(X_{t}^{A})\vert_{G}^{2} \Big] \, \textnormal{d} t
= \int_{0}^{T} I_{G}( P_{t}^{A}  \, \vert \, \mathrm{Q}) \, \textnormal{d} t < \infty.   
\]
Again, the process
\[
\mybar{M}_{s}^{A} 
= \big( \log \mybar{\ell}_{s}^{A} (\mybar{X}_{s}^{A} ) - \log \ell_{T}^{A} (X_{T}^{A} )  \big) - \mybar{F}_{s}^{A},
\qquad 0 \leqslant s \leqslant T
\]
is a martingale bounded in $L^{2}(\mathds{P})$, with representation as a stochastic integral
\[
\mybar{M}_{s}^{A}  
= \int_{0}^{s} \Big\langle \nabla \log \mybar{\ell}_{u}^{A}(\mybar{X}_{u}^{A} ) \, , \sqrt{2A(\mybar{X}_{u}^{A} )} \, \textnormal{d}\mybar{B}_{u} \Big\rangle.
\]

\subsection{The proof of \texorpdfstring{\hyperref[main.two.lis.tra]{Theorem \ref*{main.two.lis.tra}}}{Theorem 2.1}} \label{sec.2.1.lis.potmt}

We place ourselves on the filtered probability space $(\Omega,\mybar{\mathcal{F}}_{T},\mybar{\mathds{F}},\mathds{P})$. As already pointed out, our goal is to calculate the dynamics of the time-reversed relative entropy process \hyperref[eq.smdottrrop.lis]{(\ref*{eq.smdottrrop.lis})}. In order to do this, the first question is how the forward dynamics \hyperref[sde.sigma.v.lisini]{(\ref*{sde.sigma.v.lisini})} of the diffusion process $(X_{t}^{\Sigma})_{0 \leqslant t \leqslant T}$ change when regarded in the backward direction of time. In other words, we need to know the stochastic differential of the time-reversed process \hyperref[eq.trxsig.lis]{(\ref*{eq.trxsig.lis})}. This requires the classical theory of time-reversal; see, e.g., \cite{Foe85,Foe86}, \cite{HP86}, \cite{Mey94}, \cite{Nel01}, and \cite{Par86}. We summarize the results relevant for us in the following lemma. For its proof we refer to Theorems G.2 and G.5 of Appendix G in \cite{KST20}.

\begin{lemma} \label{sec.2.1.lis.potmt.lem.bm} The stochastic process $(\mybar{B}_{s})_{0 \leqslant s \leqslant T}$ given by
\[
\mybar{B}_{s} \coloneqq B_{T-s} - B_{T} - 
\int_{0}^{s} \operatorname{div}\Big( \, \mybar{p}_{u}^{\Sigma}(\mybar{X}_{u}^{\Sigma}) \, \sqrt{2\Sigma(\mybar{X}_{u}^{\Sigma})} \, \Big) \big( \, \mybar{p}_{u}^{\Sigma}(\mybar{X}_{u}^{\Sigma})\, \big)^{-1} \, \textnormal{d}u
\]
is a Brownian motion on the filtered probability space $(\Omega,\mybar{\mathcal{F}}_{T},\mybar{\mathds{F}},\mathds{P})$, and the time-reversed process \textnormal{\hyperref[eq.trxsig.lis]{(\ref*{eq.trxsig.lis})}} satisfies the stochastic differential equation
\begin{equation} \label{sde.bwd.lis.set}
\textnormal{d}\mybar{X}_{s}^{\Sigma} 
= \big( \operatorname{div}\Sigma - \Sigma  \nabla V + 2 \Sigma \, \nabla \log \mybar{\ell}_{s}^{\Sigma} \big)(\mybar{X}_{s}^{\Sigma}) \, \textnormal{d} s + \sqrt{2\Sigma(\mybar{X}_{s}^{\Sigma})} \, \textnormal{d}\mybar{B}_{s}, 
\end{equation}
for $0 \leqslant s \leqslant T$.
\end{lemma}

Here, the bar over a deterministic function means that time is reversed, i.e., the time parameter $s \in [0,T]$ is replaced by $T-s$. In the lemma above, this conventions means that $\mybar{p}_{u}^{\Sigma} = p_{T-u}^{\Sigma}$ and $\mybar{\ell}_{s}^{\Sigma} = \ell_{T-s}^{\Sigma}$.

\smallskip

The following lemma describes the partial differential equation satisfied by the likelihood ratio function \hyperref[eq.lis.lhrf.a]{(\ref*{eq.lis.lhrf.a})}. It will be crucial in order to compute the stochastic differential of the time-reversed likelihood ratio process 
\begin{equation} \label{trrep.lis.sig}
\mybar{\ell}_{s}^{\Sigma}(\mybar{X}_{s}^{\Sigma})
= \ell_{T-s}^{\Sigma} (X_{T-s}^{\Sigma} ), 
\qquad 0 \leqslant s \leqslant T
\end{equation}
and of its logarithm \hyperref[eq.smdottrrop.lis]{(\ref*{eq.smdottrrop.lis})}.

\begin{lemma} \label{lem.pde.llhrf.lis.fdot} The likelihood ratio function \textnormal{\hyperref[eq.lis.lhrf.a]{(\ref*{eq.lis.lhrf.a})}} satisfies the partial differential equation
\begin{equation} \label{lem.pde.llhrf.lis.fdot.eq.a}
\partial_{t} \ell_{t}^{\Sigma}  
= \sum_{i,j=1}^{n} \Sigma_{ij} \, \partial_{ij}^{2} \ell_{t}^{\Sigma} + \big\langle \operatorname{div} \Sigma - \Sigma \, \nabla V \, , \nabla \ell_{t}^{\Sigma}  \big\rangle.
\end{equation}
\end{lemma}

Reversing time, we derive from \hyperref[lem.pde.llhrf.lis.fdot.eq.a]{(\ref*{lem.pde.llhrf.lis.fdot.eq.a})} the equation
\begin{equation} \label{lem.pde.llhrf.lis.fdot.eq.a.dtc}
- \partial_{s} \mybar{\ell}_{s}^{\Sigma} 
= \sum_{i,j=1}^{n} \Sigma_{ij} \, \partial_{ij}^{2} \mybar{\ell}_{s}^{\Sigma}  + \big\langle \operatorname{div} \Sigma - \Sigma \,  \nabla V \, , \nabla \mybar{\ell}_{s}^{\Sigma} \big\rangle.
\end{equation}
This is exactly the \textit{Kolmogorov backward equation} corresponding to the stochastic differential equation \hyperref[sde.sigma.v.lisini]{(\ref*{sde.sigma.v.lisini})}, a circumstance indicating that the dynamics of the likelihood ratio process \hyperref[llhr.pr.itfdot.lis.s]{(\ref*{llhr.pr.itfdot.lis.s})} are most descriptive  under time reversal.

\begin{proof}[Proof of \texorpdfstring{\hyperref[lem.pde.llhrf.lis.fdot]{Lemma \ref*{lem.pde.llhrf.lis.fdot}}}] Differentiating with respect to the temporal and spatial variables the product
\[
p_{t}^{\Sigma}(x) 
= q(x) \, \ell_{t}^{\Sigma}(x) 
= \mathrm{e}^{-V(x)} \, \ell_{t}^{\Sigma}(x),
\]
substituting these derivatives into the Fokker--Planck equation \hyperref[pde.sigma.v.lisini]{(\ref*{pde.sigma.v.lisini})}, and using the stationary version \hyperref[pde.sigma.v.lisini.sv]{(\ref*{pde.sigma.v.lisini.sv})} of this equation, one obtains the partial differential equation \hyperref[lem.pde.llhrf.lis.fdot.eq.a]{(\ref*{lem.pde.llhrf.lis.fdot.eq.a})} for the likelihood ratio function.
\end{proof}

\begin{lemma} On the filtered probability space $(\Omega,\mybar{\mathcal{F}}_{T},\mybar{\mathds{F}},\mathds{P})$, the time-reversed likelihood ratio process \textnormal{\hyperref[trrep.lis.sig]{(\ref*{trrep.lis.sig})}} satisfies the stochastic differential equation
\begin{equation} \label{sde.trlrp.lis.sig}
\frac{\textnormal{d} \mybar{\ell}_{s}^{\Sigma}(\mybar{X}_{s}^{\Sigma})}{\mybar{\ell}_{s}^{\Sigma}(\mybar{X}_{s}^{\Sigma})} 
= \Big\langle  \nabla \log \mybar{\ell}_{s}^{\Sigma}(\mybar{X}_{s}^{\Sigma}) \, , \sqrt{2\Sigma(\mybar{X}_{s}^{\Sigma})} \, \textnormal{d}\mybar{B}_{s} \Big\rangle 
+ \Big\langle \nabla \log \mybar{\ell}_{s}^{\Sigma} \, , \, 2\Sigma \, \nabla \log \mybar{\ell}_{s}^{\Sigma} \Big\rangle(\mybar{X}_{s}^{\Sigma}) \, \textnormal{d}s, 
\end{equation}
for $0 \leqslant s \leqslant T$.
\begin{proof} Applying It\^{o}'s formula to the process \hyperref[trrep.lis.sig]{(\ref*{trrep.lis.sig})} and using the backward dynamics \hyperref[sde.bwd.lis.set]{(\ref*{sde.bwd.lis.set})} in conjunction with the Kolmogorov backward equation \hyperref[lem.pde.llhrf.lis.fdot.eq.a.dtc]{(\ref*{lem.pde.llhrf.lis.fdot.eq.a.dtc})}, we obtain the stochastic differential equation \hyperref[sde.trlrp.lis.sig]{(\ref*{sde.trlrp.lis.sig})}.
\end{proof}
\end{lemma}

\begin{lemma} \label{lem.lis.log.retr.sig} On the filtered probability space $(\Omega,\mybar{\mathcal{F}}_{T},\mybar{\mathds{F}},\mathds{P})$, the time-reversed relative entropy process \textnormal{\hyperref[eq.smdottrrop.lis]{(\ref*{eq.smdottrrop.lis})}} satisfies the stochastic differential equation
\begin{equation} \label{sde.trrep.lis.sig}
\textnormal{d} \log \mybar{\ell}_{s}^{\Sigma}(\mybar{X}_{s}^{\Sigma}) 
= \Big\langle  \nabla \log \mybar{\ell}_{s}^{\Sigma}(\mybar{X}_{s}^{\Sigma}) \, , \sqrt{2\Sigma(\mybar{X}_{s}^{\Sigma})} \, \textnormal{d}\mybar{B}_{s} \Big\rangle 
+ \Big\langle \nabla \log \mybar{\ell}_{s}^{\Sigma} \, , \, \Sigma \, \nabla \log \mybar{\ell}_{s}^{\Sigma} \Big\rangle(\mybar{X}_{s}^{\Sigma}) \, \textnormal{d}s, 
\end{equation}
for $0 \leqslant s \leqslant T$.
\begin{proof} We just have to apply It\^{o}'s formula to the process \hyperref[eq.smdottrrop.lis]{(\ref*{eq.smdottrrop.lis})} and use the dynamics \hyperref[sde.trlrp.lis.sig]{(\ref*{sde.trlrp.lis.sig})}.
\end{proof}
\end{lemma}

With \hyperref[lem.lis.log.retr.sig]{Lemma \ref*{lem.lis.log.retr.sig}} we have prepared the main ingredient for the proof of \hyperref[main.two.lis.tra]{Theorem \ref*{main.two.lis.tra}}. In order to prove the integrability condition \hyperref[main.two.lis.tra.a]{(\ref*{main.two.lis.tra.a})} we follow a similar strategy as in the proof of Theorem 4.1 in \cite{KST22}.

\begin{proof}[\bfseries \upshape Proof of \texorpdfstring{\hyperref[main.two.lis.tra]{Theorem \ref*{main.two.lis.tra}}}{}] The fact that the process $\mybar{M}^{\Sigma} = (\mybar{M}_{s}^{\Sigma})_{0 \leqslant s \leqslant T}$ of \hyperref[main.two.lis.tra.b]{(\ref*{main.two.lis.tra.b})} is a local martingale with stochastic integral representation \hyperref[main.two.lis.tra.c]{(\ref*{main.two.lis.tra.c})} follows from the stochastic differential \hyperref[sde.trrep.lis.sig]{(\ref*{sde.trrep.lis.sig})} established in \hyperref[lem.lis.log.retr.sig]{Lemma \ref*{lem.lis.log.retr.sig}}. It remains to show the integrability condition \hyperref[main.two.lis.tra.a]{(\ref*{main.two.lis.tra.a})}, which in particular implies that $\mybar{M}^{\Sigma}$ is a martingale bounded in $L^{2}(\mathds{P})$.

\smallskip

To this end, let us consider the inverse of the process \hyperref[trrep.lis.sig]{(\ref*{trrep.lis.sig})}, i.e.,
\begin{equation} \label{llhr.pr.itfdot.lis.s.inv}
(\mybar{\ell}_{s}^{\Sigma})^{-1}(\mybar{X}_{s}^{\Sigma})
= \frac{q(\mybar{X}_{s}^{\Sigma})}{\mybar{p}_{s}^{\Sigma}(\mybar{X}_{s}^{\Sigma})}, 
\qquad 0 \leqslant s \leqslant T.
\end{equation}
Applying It\^{o}'s formula and using \hyperref[sde.trlrp.lis.sig]{(\ref*{sde.trlrp.lis.sig})}, we obtain the stochastic differential
\begin{equation} \label{sde.trlrp.lis.sig.inv.d}
\textnormal{d} (\mybar{\ell}_{s}^{\Sigma})^{-1}(\mybar{X}_{s}^{\Sigma})
= - (\mybar{\ell}_{s}^{\Sigma})^{-1}(\mybar{X}_{s}^{\Sigma}) \,
\Big\langle  \nabla \log \mybar{\ell}_{s}^{\Sigma}(\mybar{X}_{s}^{\Sigma}) \, , \sqrt{2\Sigma(\mybar{X}_{s}^{\Sigma})} \, \textnormal{d}\mybar{B}_{s} \Big\rangle 
\end{equation}
and conclude that the process \hyperref[llhr.pr.itfdot.lis.s.inv]{(\ref*{llhr.pr.itfdot.lis.s.inv})} is a local martingale. From \hyperref[sde.trlrp.lis.sig.inv.d]{(\ref*{sde.trlrp.lis.sig.inv.d})}, or directly from \hyperref[sde.trrep.lis.sig]{(\ref*{sde.trrep.lis.sig})}, we obtain the dynamics of the logarithm of the process \hyperref[llhr.pr.itfdot.lis.s.inv]{(\ref*{llhr.pr.itfdot.lis.s.inv})}, namely
\[
\textnormal{d} \log (\mybar{\ell}_{s}^{\Sigma})^{-1}(\mybar{X}_{s}^{\Sigma}) 
= - \Big\langle  \nabla \log \mybar{\ell}_{s}^{\Sigma}(\mybar{X}_{s}^{\Sigma}) \, , \sqrt{2\Sigma(\mybar{X}_{s}^{\Sigma})} \, \textnormal{d}\mybar{B}_{s} \Big\rangle 
-  \Big\langle \nabla \log \mybar{\ell}_{s}^{\Sigma} \, , \, \Sigma \, \nabla \log \mybar{\ell}_{s}^{\Sigma} \Big\rangle(\mybar{X}_{s}^{\Sigma}) \, \textnormal{d}s.
\]
At this point, let us recall the relative entropy \hyperref[re.kld.lisini]{(\ref*{re.kld.lisini})} and point \hyperref[def.cl.of.ad.vol.sig.lis.vi]{\ref*{def.cl.of.ad.vol.sig.lis.vi}} of \hyperref[def.cl.of.ad.vol.sig.lis]{Definition \ref*{def.cl.of.ad.vol.sig.lis}}, which ensures that the initial relative entropy $H( P_{0}^{\Sigma} \, \vert \, \mathrm{Q} )$ is finite. This implies that the terminal value $\log (\mybar{\ell}_{T}^{\Sigma})^{-1}(\mybar{X}_{T}^{\Sigma})$ is integrable, i.e.,
\[
\mathds{E}_{\mathds{P}}\big[ \log (\mybar{\ell}_{T}^{\Sigma})^{-1}(\mybar{X}_{T}^{\Sigma})\big] 
= - H( P_{0}^{\Sigma} \, \vert \, \mathrm{Q} ) \in (-\infty,\infty).
\]
On the other hand, for the initial value we have that 
\begin{equation} \label{eq.lis.iv.f.yt}
\mathds{E}_{\mathds{P}}\big[ \log (\mybar{\ell}_{0}^{\Sigma})^{-1}(\mybar{X}_{0}^{\Sigma})\big] 
= - H( P_{T}^{\Sigma} \, \vert \, \mathrm{Q} ) \in [-\infty,\infty).
\end{equation}
Thus we can apply \cite[Proposition A.3]{KST22} to the local martingale \hyperref[llhr.pr.itfdot.lis.s.inv]{(\ref*{llhr.pr.itfdot.lis.s.inv})} and the deterministic stopping time $\tau = T$. We conclude that the expectation in \hyperref[eq.lis.iv.f.yt]{(\ref*{eq.lis.iv.f.yt})} is finite and we have the identity
\[
\mathds{E}_{\mathds{P}}\big[ \log (\mybar{\ell}_{T}^{\Sigma})^{-1}(\mybar{X}_{T}^{\Sigma})\big] 
- \mathds{E}_{\mathds{P}}\big[ \log (\mybar{\ell}_{0}^{\Sigma})^{-1}(\mybar{X}_{0}^{\Sigma})\big] 
= - \mathds{E}_{\mathds{P}} \bigg[ \int_{0}^{T} \Big\langle \nabla \log \mybar{\ell}_{u}^{\Sigma} \, , \, \Sigma \, \nabla \log \mybar{\ell}_{u}^{\Sigma} \Big\rangle(\mybar{X}_{u}^{\Sigma}) \, \textnormal{d}u \bigg].
\]
This implies that the local martingale $\mybar{M}^{\Sigma}$ is bounded in $L^{2}(\mathds{P})$. In fact, by the It\^{o} isometry, the $L^{2}(\mathds{P})$-norm of the stochastic integral
\[
\mybar{M}_{T}^{\Sigma}  
= \int_{0}^{T} \Big\langle \nabla \log \mybar{\ell}_{u}^{\Sigma}(\mybar{X}_{u}^{\Sigma} ) \, , \sqrt{2\Sigma(\mybar{X}_{u}^{\Sigma} )} \, \textnormal{d}\mybar{B}_{u} \Big\rangle
\]
is finite and can be computed as
\[
\tfrac{1}{2} \big\Vert \mybar{M}_{T}^{\Sigma} \big\Vert_{L^{2}(\mathds{P})}^{2}
= \mathds{E}_{\mathds{P}} \bigg[ \int_{0}^{T} \Big\langle \nabla \log \ell_{t}^{\Sigma} \, , \, \Sigma \, \nabla \log \ell_{t}^{\Sigma} \Big\rangle(X_{t}^{\Sigma}) \, \textnormal{d}t \bigg]
= H( P_{0}^{\Sigma} \, \vert \, \mathrm{Q} )
- H( P_{T}^{\Sigma} \, \vert \, \mathrm{Q} ).
\]
Recalling \hyperref[cum.rel.fis.inf.pro.lis.sig]{(\ref*{cum.rel.fis.inf.pro.lis.sig})}, this shows the integrability condition \hyperref[main.two.lis.tra.a]{(\ref*{main.two.lis.tra.a})} and finishes the proof of \hyperref[main.two.lis.tra]{Theorem \ref*{main.two.lis.tra}}.
\end{proof}

\subsection{Consequences of \texorpdfstring{\hyperref[main.two.lis.tra]{Theorem \ref*{main.two.lis.tra}}}{Theorem 2.1}} \label{cottt.lis.ram.sig.a}

In this subsection we state several important consequences of our basic trajectorial result, \hyperref[main.two.lis.tra]{Theorem \ref*{main.two.lis.tra}}.

\begin{corollary} \label{cottt.lis.ram.sig.a.cor.one} Let $\Sigma \in \mathcal{S}$, fix $A \in \mathcal{S}$ and set $G = A^{-1}$. For all $0 \leqslant t_{0},t \leqslant T$ we have the relative entropy identity
\begin{equation} \label{eq.cottt.lis.ram.sig.a.cor.one}
H( P_{t_{0}}^{\Sigma} \, \vert \, \mathrm{Q}) - H( P_{t}^{\Sigma} \, \vert \, \mathrm{Q}) 
= \int_{t_{0}}^{t} \mathds{E}_{\mathds{P}}\Big[\Big\langle \nabla \log \ell_{u}^{\Sigma} \, , \, \Sigma  \, \nabla \log \ell_{u}^{\Sigma} \Big\rangle(X_{u}^{\Sigma}) \Big] \, \textnormal{d}u.
\end{equation}
Furthermore, for Lebesgue-a.e.\ $0 \leqslant t \leqslant T$, the dissipation of relative entropy is equal to
\begin{equation} \label{eq.cottt.lis.ram.sig.a.cor.two}
\frac{\textnormal{d}}{\textnormal{d}t} \,  H( P_{t}^{\Sigma} \, \vert \, \mathrm{Q}) 
= \mathds{E}_{\mathds{P}}\Big[\Big\langle \nabla \log \ell_{t}^{\Sigma} \, , \, \Sigma \, \nabla \log \ell_{t}^{\Sigma} \Big\rangle(X_{t}^{\Sigma}) \Big].
\end{equation}
\begin{proof} The identity \hyperref[eq.cottt.lis.ram.sig.a.cor.one]{(\ref*{eq.cottt.lis.ram.sig.a.cor.one})} follows directly from \hyperref[main.two.lis.tra]{Theorem \ref*{main.two.lis.tra}} by taking expectations and using that the process $\mybar{M}^{\Sigma}$ of \hyperref[main.two.lis.tra.b]{(\ref*{main.two.lis.tra.b})} is a martingale. Applying the Lebesgue differentiation theorem yields \hyperref[eq.cottt.lis.ram.sig.a.cor.two]{(\ref*{eq.cottt.lis.ram.sig.a.cor.two})}.
\end{proof}
\end{corollary}

As a consequence of \hyperref[eq.cottt.lis.ram.sig.a.cor.one]{(\ref*{eq.cottt.lis.ram.sig.a.cor.one})}, we observe that the relative entropy function (recall \hyperref[rel.ent.function.lis.sig]{(\ref*{rel.ent.function.lis.sig})})
\[
[0,T] \ni t \longmapsto H( P_{t}^{\Sigma} \, \vert \, \mathrm{Q})
\]
is non-increasing. As the initial relative entropy $H( P_{0}^{\Sigma} \, \vert \, \mathrm{Q})$ is finite by \hyperref[def.cl.of.ad.vol.sig.lis]{Definition \ref*{def.cl.of.ad.vol.sig.lis}}, \hyperref[def.cl.of.ad.vol.sig.lis.vi]{\ref*{def.cl.of.ad.vol.sig.lis.vi}}, this implies that the relative entropy $H( P_{t}^{\Sigma} \, \vert \, \mathrm{Q})$ is finite, for all times $0 \leqslant t \leqslant T$.

\smallskip

Choosing $\Sigma = A$ and recalling the relative Fisher information of \hyperref[rel.fis.inf.det.pexp.sig]{(\ref*{rel.fis.inf.det.pexp.sig})}, we immediately obtain from \hyperref[cottt.lis.ram.sig.a.cor.one]{Corollary \ref*{cottt.lis.ram.sig.a.cor.one}} the following result.

\begin{corollary} \label{cottt.lis.ram.A.a.cor.one} Let $A \in \mathcal{S}$ and set $G = A^{-1}$. For all $0 \leqslant t_{0},t \leqslant T$ we have the relative entropy identity
\begin{equation} \label{eq.cottt.lis.ram.A.a.cor.one}
H( P_{t_{0}}^{A} \, \vert \, \mathrm{Q}) - H( P_{t}^{A} \, \vert \, \mathrm{Q}) 
= 
\int_{t_{0}}^{t} \mathds{E}_{\mathds{P}}\Big[\Big\langle \nabla \log \ell_{u}^{A} \, , \, A  \, \nabla \log \ell_{u}^{A} \Big\rangle(X_{u}^{A}) \Big] \, \textnormal{d}u 
= \int_{t_{0}}^{t} I_{G}( P_{u}^{A} \, \vert \, \mathrm{Q})  \, \textnormal{d}u.
\end{equation}
Furthermore, for Lebesgue-a.e.\ $0 \leqslant t \leqslant T$, the dissipation of relative entropy is equal to
\begin{equation} \label{eq.cottt.lis.ram.A.a.cor.two}
\frac{\textnormal{d}}{\textnormal{d}t} \,  H( P_{t}^{A} \, \vert \, \mathrm{Q}) 
= \mathds{E}_{\mathds{P}}\Big[\Big\langle \nabla \log \ell_{t}^{A} \, , \, A  \, \nabla \log \ell_{t}^{A} \Big\rangle(X_{t}^{A}) \Big] \\
=  I_{G}( P_{t}^{A} \, \vert \, \mathrm{Q}).
\end{equation}
\end{corollary}

The \textit{entropy dissipation} \hyperref[eq.cottt.lis.ram.A.a.cor.two]{(\ref*{eq.cottt.lis.ram.A.a.cor.two})} is also is also known as an \textit{entropy production equality}. Various authors have studied this identity in different settings, let us refer to \cite{BE85, OV00, AMTU01, CMV03, Vil03, AGS08, Lis09, KST22}. The relation 
\[
\frac{\textnormal{d}}{\textnormal{d}t} \,  H( P_{t}^{A} \, \vert \, \mathrm{Q}) 
=  I_{G}( P_{t}^{A} \, \vert \, \mathrm{Q})
\]
between entropy dissipation and Fisher information also takes the familiar form of a \textit{de Bruijn type identity} from information theory (see, e.g. \cite{CT06}).

\smallskip

We now return to the probabilistic setting and place ourselves on the filtered probability space $(\Omega,\mybar{\mathcal{F}}_{T},\mybar{\mathds{F}},\mathds{P})$ as in \hyperref[main.two.lis.tra]{Theorem \ref*{main.two.lis.tra}}. The following result provides the trajectorial analogue of \hyperref[eq.cottt.lis.ram.sig.a.cor.two]{(\ref*{eq.cottt.lis.ram.sig.a.cor.two})}.

\begin{corollary} Let $\Sigma \in \mathcal{S}$, fix $A \in \mathcal{S}$ and set $G = A^{-1}$. For Lebesgue-a.e.\ $0 \leqslant t \leqslant T$, the trajectorial rate of relative entropy dissipation is equal to
\begingroup
\addtolength{\jot}{0.7em}
\begin{align} 
\lim_{s \uparrow T-t} 
\frac{\log \mybar{\ell}_{s}^{\Sigma}(\mybar{X}_{s}^{\Sigma}) - \mathds{E}_{\mathds{P}}\big[ \log \ell_{t}^{\Sigma}(X_{t}^{\Sigma}) \, \vert \, \mybar{\mathcal{F}}_{s} \big]}{T-t-s}  
&= - \Big\langle \nabla \log \ell_{t}^{\Sigma} \, , \, \Sigma  \, \nabla \log \ell_{t}^{\Sigma} \Big\rangle(X_{t}^{\Sigma}), \label{trored.eq.cor.lis.sig.a} \\
\lim_{s \downarrow T-t} 
\frac{\log \ell_{t}^{\Sigma}(X_{t}^{\Sigma})-\mathds{E}_{\mathds{P}}\big[ \log \mybar{\ell}_{s}^{\Sigma}(\mybar{X}_{s}^{\Sigma}) \, \vert \, \mybar{\mathcal{F}}_{T-t} \big]}{s-(T-t)}   
&= - \Big\langle \nabla \log \ell_{t}^{\Sigma} \, , \, \Sigma  \, \nabla \log \ell_{t}^{\Sigma} \Big\rangle(X_{t}^{\Sigma}), \label{trored.eq.cor.lis.sig.b}
\end{align}
\endgroup
where both limits exist in $L^{1}(\mathds{P})$.
\begin{proof} It follows from \hyperref[main.two.lis.tra]{Theorem \ref*{main.two.lis.tra}} that the numerator of the fraction on the left-hand side of \hyperref[trored.eq.cor.lis.sig.a]{(\ref*{trored.eq.cor.lis.sig.a})} is equal to 
\begin{equation} \label{trored.eq.cor.lis.sig.a.i}
- \mathds{E}_{\mathds{P}} \Big[ \mybar{F}_{T-t}^{\Sigma} - \mybar{F}_{s}^{\Sigma} \, \big\vert \, \mybar{\mathcal{F}}_{s} \Big] 
= - \mathds{E}_{\mathds{P}}\bigg[\int_{s}^{T-t} \Big\langle \nabla \log \mybar{\ell}_{u}^{\Sigma} \, , \, \Sigma \, \nabla \log \mybar{\ell}_{u}^{\Sigma} \Big\rangle(\mybar{X}_{u}^{\Sigma}) \, \textnormal{d}u \ \Big\vert \ \mybar{\mathcal{F}}_{s} \bigg].
\end{equation}
After dividing by $T-t-s$, the $L^{1}(\mathds{P})$-convergence to the right-hand side of \hyperref[trored.eq.cor.lis.sig.a]{(\ref*{trored.eq.cor.lis.sig.a})}, as $s \uparrow T - t$, follows from the Lebesgue differentiation theorem and its generalization to conditional expectations (see, e.g., \cite[Proposition A.2]{KST22}). The proof of \hyperref[trored.eq.cor.lis.sig.b]{(\ref*{trored.eq.cor.lis.sig.b})} follows a similar line of reasoning.
\end{proof}
\end{corollary}

\section{Proof of the gradient flow property} \label{ch.lis.tpotgrp}

In this section we prove the gradient flow property as formulated in \hyperref[main.one.lis.det]{Theorem \ref*{main.one.lis.det}}. In a nutshell, we show that the flow of probability measures $(P_{t}^{A})_{0 \leqslant t \leqslant T}$, defined by the stochastic differential equation \hyperref[sde.sigma.v.lisini.A]{(\ref*{sde.sigma.v.lisini.A})}, is the gradient flow of the relative entropy functional \hyperref[eq.rel.ent.func.lis.]{(\ref*{eq.rel.ent.func.lis.})} with respect to the quadratic Wasserstein distance $W_{2,G}$ as defined in \hyperref[def.eq.wtwog.lis.was]{(\ref*{def.eq.wtwog.lis.was})}.

\subsection{Derivative of the Wasserstein distance} \label{sec.3.lis.dotwd}

We define the vector field
\[
[0,T] \times \mathds{R}^{n} \longrightarrow \mathds{R}^{n} \colon (t,x) \longmapsto 
\boldsymbol{v}_{t}^{\Sigma}(x) \coloneqq - \Sigma(x) \, \big( \nabla \log p_{t}^{\Sigma}(x) + \nabla V(x)\big).
\]
Then we can write the Fokker--Planck equation \hyperref[pde.sigma.v.lisini]{(\ref*{pde.sigma.v.lisini})} as a \textit{continuity equation}, to wit
\[
\partial_{t} p_{t}^{\Sigma}(x) 
+ \operatorname{div}\big(  \boldsymbol{v}_{t}^{\Sigma}(x) \, p_{t}^{\Sigma}(x) \big)
= 0,
\qquad (t,x) \in  (0,T) \times \mathds{R}^{n}.
\]
The derivative of the Wasserstein distance along an absolutely continuous ``curve of probability measures'' satisfying a continuity equation is a well known result; see, e.g.\ \cite{AGS08, Lis09, San15}.

\begin{lemma} \label{lem.lis.ags.dotwd.a} Let $\Sigma \in \mathcal{S}$, fix $A \in \mathcal{S}$ and set $G = A^{-1}$. For Lebesgue-a.e.\ $t \in [0,T]$, we have that
\begin{equation} \label{eq.lem.lis.ags.dotwd.a}
\lim_{h \rightarrow 0} \frac{W_{2,G}(P_{t+h}^{\Sigma},P_{t}^{\Sigma})}{\vert h \vert} 
= \Vert \boldsymbol{v}_{t}^{\Sigma} \Vert_{L_{G}^{2}(P_{t}^{\Sigma})}
= \sqrt{\mathds{E}_{\mathds{P}}\Big[ \big\vert \big(\Sigma \, G \, \nabla_{G} \log \ell_{t}^{\Sigma}\big)(X_{t}^{\Sigma})\big\vert_{G}^{2}  \Big]}.
\end{equation}
\begin{proof} For a proof of this result we refer to \cite[Theorem 2.4]{Lis09}, which is an extension of \cite[Theorem 8.3.1]{AGS08} to the current Riemannian setting induced by $G$. Alternatively, the reader may consult \cite[Theorem 5.14]{San15} and adapt the arguments by means of \cite{McC01}, so that they apply also in our Riemannian setting.
\end{proof}
\end{lemma}

In the case $\Sigma = A$, the derivative \hyperref[eq.lem.lis.ags.dotwd.a]{(\ref*{eq.lem.lis.ags.dotwd.a})} takes the form
\begin{equation} \label{eq.lem.lis.ags.dotwd.a.sig.is.a}
\lim_{h \rightarrow 0} \frac{W_{2,G}(P_{t+h}^{A},P_{t}^{A})}{\vert h \vert} 
= \Vert \boldsymbol{v}_{t}^{A} \Vert_{L_{G}^{2}(P_{t}^{A})}
= \sqrt{\mathds{E}_{\mathds{P}}\big[ \vert \nabla_{G} \log \ell_{t}^{A}(X_{t}^{A})\vert_{G}^{2}  \big]}
= \sqrt{I_{G}( P_{t}^{A} \, \vert \, \mathrm{Q})};
\end{equation}
for the last equality in \hyperref[eq.lem.lis.ags.dotwd.a.sig.is.a]{(\ref*{eq.lem.lis.ags.dotwd.a.sig.is.a})}, recall the definition of relative Fisher information \hyperref[rel.fis.inf.det.pexp.sig]{(\ref*{rel.fis.inf.det.pexp.sig})}.

\subsection{The proof of \texorpdfstring{\hyperref[main.one.lis.det]{Theorem \ref*{main.one.lis.det}}}{Theorem 1.2}} \label{section4.lis.tpot11}

\begin{proof}[\bfseries \upshape Proof of \texorpdfstring{\hyperref[main.one.lis.det]{Theorem \ref*{main.one.lis.det}}}{}] Let $\Sigma \in \mathcal{S}$, fix $A \in \mathcal{S}$ and set $G = A^{-1}$. Let $0 \leqslant t_{0} \leqslant t \leqslant T$. We write the relative entropy identity
\hyperref[eq.cottt.lis.ram.sig.a.cor.one]{(\ref*{eq.cottt.lis.ram.sig.a.cor.one})} of \hyperref[cottt.lis.ram.sig.a.cor.one]{Corollary \ref*{cottt.lis.ram.sig.a.cor.one}} as 
\[
H( P_{t_{0}}^{\Sigma} \, \vert \, \mathrm{Q}) - H( P_{t}^{\Sigma} \, \vert \, \mathrm{Q}) 
= \int_{t_{0}}^{t} \mathds{E}_{\mathds{P}}\Big[\Big\langle \nabla_{G} \log \ell_{u}^{\Sigma} \, , \, \Sigma  \, G \, \nabla_{G} \log \ell_{u}^{\Sigma} \Big\rangle_{G}(X_{u}^{\Sigma}) \Big] \, \textnormal{d}u. 
\]
Applying the Cauchy--Schwarz inequality yields
\[
\big\vert H( P_{t}^{\Sigma} \, \vert \, \mathrm{Q}) - H( P_{t_{0}}^{\Sigma} \, \vert \, \mathrm{Q})\big\vert
\leqslant 
\int_{t_{0}}^{t} \sqrt{\mathds{E}_{\mathds{P}}\Big[ \big\vert \nabla_{G} \log \ell_{u}^{\Sigma}(X_{u}^{\Sigma})\big\vert_{G}^{2}  \Big]} \
\sqrt{\mathds{E}_{\mathds{P}}\Big[ \big\vert \big(\Sigma \, G \, \nabla_{G} \log \ell_{u}^{\Sigma}\big)(X_{u}^{\Sigma})\big\vert_{G}^{2}  \Big]} \, \textnormal{d}u.
\]
Recalling the definition of relative Fisher information \hyperref[rel.fis.inf.det.pexp.sig]{(\ref*{rel.fis.inf.det.pexp.sig})} and the derivative of the Wasserstein distance \hyperref[eq.lem.lis.ags.dotwd.a]{(\ref*{eq.lem.lis.ags.dotwd.a})}, this proves the inequality \hyperref[main.one.lis.det.01]{(\ref*{main.one.lis.det.01})}.
Clearly, if $\Sigma$ is a constant multiple of $A$, this inequality is an equality.

\smallskip

Now let us consider the case $\Sigma = A$. Recalling the relative entropy identity \hyperref[eq.cottt.lis.ram.A.a.cor.one]{(\ref*{eq.cottt.lis.ram.A.a.cor.one})} of \hyperref[cottt.lis.ram.A.a.cor.one]{Corollary \ref*{cottt.lis.ram.A.a.cor.one}} we have
\[
H( P_{t_{0}}^{A} \, \vert \, \mathrm{Q}) - H( P_{t}^{A} \, \vert \, \mathrm{Q}) 
= \int_{t_{0}}^{t} I_{G}( P_{u}^{A} \, \vert \, \mathrm{Q})  \, \textnormal{d}u.
\]
According to \hyperref[eq.lem.lis.ags.dotwd.a.sig.is.a]{(\ref*{eq.lem.lis.ags.dotwd.a.sig.is.a})}, we have the derivative of the Wasserstein distance 
\[
\lim_{h \rightarrow 0} \frac{W_{2,G}(P_{t+h}^{A},P_{t}^{A})}{\vert h \vert} 
= \sqrt{I_{G}( P_{t}^{A} \, \vert \, \mathrm{Q})},
\]
so that we can write
\[
I_{G}( P_{u}^{A} \, \vert \, \mathrm{Q}) = 
\tfrac{1}{2} I_{G}( P_{u}^{A} \, \vert \, \mathrm{Q}) + \tfrac{1}{2} \bigg( \lim_{h \rightarrow 0} \frac{W_{2,G}(P_{u+h}^{A},P_{u}^{A})}{\vert h \vert} \bigg)^{2}.
\]
We conclude the equality \hyperref[main.one.lis.det.02]{(\ref*{main.one.lis.det.02})}.
\end{proof}

\bibliographystyle{alpha}
{\footnotesize
\bibliography{references}}

\end{document}